\documentclass[12pt]{article}
\usepackage{graphicx}
\usepackage{amssymb}

\begin{document}

\title{An ecosystem with Holling type II response and predators' genetic variability
}

\author{Clara Viberti,
Ezio Venturino\thanks{Corresponding author.
Email: ezio.venturino@unito.it}\\
Dipartimento di Matematica ``Giuseppe Peano'',\\
Universit\`{a} di Torino, Italy.
}
\date{}
\maketitle

\begin{abstract}
A new model to investigate environmental effects of genetically distinguishable predators is presented.
The Holling type II response function, modelling feeding satiation, leads to persistent system's oscillations,
as in classical population models. An almost complete classification of the cases arising in the Routh-Hurwitz stability
conditions mathematically characterizes the paper. It is instrumental as a guideline in the numerical experiments leading to
the findings on the limit cycles. This result extends what found in an earlier parallel investigation containing a standard bilinear
response function.
\end{abstract}



\section{Introduction}

Mathematical models for the interactions of populations are now classical, \cite{Mu,BCC}.
Studies of interacting populations among which diseases spread constitute the object of ecoepidemiology, which dates back to a paper
of the late 1980's, \cite{HF}, and progressed through early works in different biological settings, such as predator-prey models,
\cite{EV92,EV94,EV95,CA,CB,EV02}, oceanic environment, \cite{BC,CP}, competing and symbiotic interactions \cite{EV01,EV07}.
The interested reader can consult Chapter 7 of \cite{MPV} for a
fairly recent account on progress of this discipline.

An extension of this situation has recently been proposed, in which the disease is
not a basic fundamental ingredient of the ecosystem, but it is replaced by the presence of more than one genotype in one of the populations.
In a sense, however, epidemics keep on playing a role in this context, since the genotype itself may make the indivuals carrying it more prone to
a certain specific disease. In this respect, these systems are very much related to ecoepidemiology. They could be referred to as mathematical
ecogenetics models, since ecogenetics is very well established discipline of biology.
Indeed, it mainly investigates how an inherited genetical variability responds to environmental changes, such as substances present in it, \cite{Cal,Cu}.
Our focus lies instead on the ecosystem behavioral consequences of the presence and interplay of
the genetically interesting population with the other populations.

The case of a genetically differentiated prey population subject to predation by their
natural predators has been presented and analysed in \cite{EV12}. Models for natural situations in which more predators feed on the
same prey are well known, \cite{LC07}.
In \cite{VV} therefore the study has been extended to the case in which the predators show genetic differences. Bilinear interaction terms have
been assumed, corresponding to the standard quadratic model in population theory. No population oscillations have been discovered.
Here we continue the investigation, in the search for possible interesting features in the system behavior.
The model is thus formulated using a Holling type II response function, as the latter better suited to model feeding, which is subject to satiation
when too large amounts of prey are present, \cite{LC03,LC06}.
In view of the large number of parameters of the model, a blind search in the parameter space for
a configuration that leads to persistent oscillations is very difficult. However, we do provide
an almost complete classification of all the cases that can arise. This mathematical effort specifically characterizes this
investigation. It is instrumental to provide guidelines for the parameter choices, and its usefulness is shown by the fact that
on this basis limit cycles are indeed found in the numerical simulations.

The paper is organized as follows. The next Section contains the model. The equilibria are studied in Section 3. A thorough classification of
the Routh-Hurwitz conditions in terms of the model parameters is carried out in Section 4. The following Section contains numerical examples
that have been worked out on the experience matured in constructing the previous classification. A brief final discussion concludes the paper.

\section{The model}

Let the predators be genetically diversified, with the two genotypes denoted by $y(\tau)$ and $z(\tau)$, let $x(\tau)$ be the prey population.
We consider the following model
\begin{eqnarray}\label{H1}
\displaystyle
\ \ x'(\tau)=R\left(1-\frac{x}{\widetilde K}\right)x-h\displaystyle\frac{\xi  x}{x+\mu}y-g\displaystyle \frac{\xi  x}{x+\mu} z,\\
\displaystyle \nonumber
\ \ y'(\tau)=pe(hy+gz)\displaystyle\frac{\xi  x}{x+\mu}-my, \\
\displaystyle \nonumber
\ \ z'(\tau)=qe(hy+gz)\displaystyle \frac{\xi  x}{x+\mu}-nz.
\end{eqnarray}
Here all the parameters are always assumed to be nonnegative.

In this situation, the key factor is here represented by the term in bracket in the last two equations.
It contains both genotypes, meaning that it is the whole predator population that reproduces. But furthermore,
since both subpopulations appear as reproduction factors in both predators' equations,
this term states that each genotype can give rise to newborns of both genotypes, where
$p$ and $q$ denote the fractions of $y$ and $z$ newborn predators, with $p+q=1$.
In other words,
the fundamental point of this model, that singles it out from other standard models in population theory,
states that from the subpopulation $y$, newborns of genotype $z$
can be generated, and vice versa. The reason for the appearance of the term $hy+gz$ can however
be pointed out more precisely as follows.
For more generality, the two genotypes are assumed to possibly have different hunting capabilities, here represented
by the coefficients $h$ for $y$ and $g$ for $z$. Each predator subpopulation thus independently
removes prey at its own rate. We also assume that the
predators experience a feeding saturation effect,
which is suitably modeled by a Holling type II response function, where
$\xi$ represents the maximum obtainable resource from each prey per unit time and
$\mu$ denotes the half saturation constant.
The total benefit from hunting for the predator population is thus represented
by the sum of these separate removing contributions for the two subpopulations,
therefore giving rise to the first term in the last two equations.
Further, newborns are produced
by converting the captured prey into new predator biomass, $e<1$ being the conversion factor.

The remaining assumptions are kind of standard in interacting population models. Namely,
the predators dynamics further shows a natural mortality, at rates $m$ and $n$ respectively for $y$ and $z$.
The prey reproduce logistically with rate $R$ and carrying capacity $\widetilde K$
and are subject to hunting by the predators, as explained above.

The system (\ref{H1}) can be nondimensionalized in the following way. Let
$\displaystyle x(\tau)=\alpha X(t)$, $\displaystyle y(\tau)=\beta Y(t)$, $\displaystyle z(\tau)=\gamma Z(t)$ e $\displaystyle t=\delta \tau$,
and choosing $\alpha=\widetilde K$, $\beta=\gamma=\displaystyle\frac{e}{g}$, $\delta=e$, we can define the new parameters
$$
r=\displaystyle\frac{R}{e}, \quad c=\displaystyle\frac{h}{g}, \quad w=pg\widetilde K, \quad s=\displaystyle\frac{m}{e},
\quad v=qg\widetilde K, \quad d=\displaystyle\frac{n}{e}.
$$
Finally, by letting
$B=\xi \widetilde K^{-1}$ and $A=\mu \widetilde K^{-1}$, we have the rescaled model
\begin{eqnarray}\label{H2}
\displaystyle
\ \ X'(t)=r (1-X)X- c\displaystyle \frac{B X}{X+A}Y - \displaystyle \frac{B X}{X+A} Z\\ \nonumber
\displaystyle
\ \ Y'(t)=w (cY+Z)\displaystyle \frac{B X}{X+A}- sY \\ \nonumber
\displaystyle
\ \ Z'(t)=v(cY+Z)\displaystyle \frac{B X}{X+A} -dZ
\end{eqnarray}

\section[Punti di equilibrio]{Equilibria}
The model (\ref{H2}) has only three possible equilibria, the origin $F_0$,
corresponding to the system extinction, the predator-free equilibrium $F_1=(1,0,0)$ and the whole ecosystem coexistence $F_2=(X^*, Y^*, Z^*)$,
the population levels of which are obtained solving for $Z$ from the first equation, substituting it into the second one to give $Y$, with
the final substitution into the last equation. This last step produces a factored quadratic, from which once again the equilibrium $F_1$ is found,
or alternatively by back substitution, the following values for the coexisting populations are determined,
$$
X^*=\displaystyle \frac{Ads}{V}, \quad
Y^*=\displaystyle \frac{wAdrW}{V^2}, \quad
Z^*=\displaystyle\frac{vAsrW}{V^2}.
$$
with 
$$
V= BQ-ds, \quad W=BQ-ds(A+1), \quad Q=sv+cdw.
$$
The feasibility conditions for $F_2$ are $B(sv+cdw)>ds$, i.e. $V>0$, and $ds(A+1)\le B(sv+cdw)$, i.e. $W\ge 0$, which combine to give
\begin{equation}\label{sistH1}
B(sv+cdw)\ge ds(A+1) \equiv \max\left\{ds,ds(A+1)\right\}.
\end{equation}

The Jacobian of (\ref{H2}) reads
\begin{equation}
\label{JH}
J=\displaystyle
\left[\begin{array}
{ccc}
r(1-2X)-\displaystyle\frac{(cY+Z)B}{A+X}+\displaystyle \frac{(cY+Z)BX}{(A+X)^2}& -\displaystyle\frac{cBX}{A+X}& -\displaystyle\frac{BX}{A+X}\\ \\
\displaystyle\frac{w(cY+Z)B}{A+X}-\displaystyle\frac{w(cY+Z)BX}{(A+X)^2} & \displaystyle\frac{wcBX}{A+X}-s & \displaystyle\frac{wBX}{A+X} \\ \\
\displaystyle\frac{v(cY+Z)B}{A+X}-\displaystyle\frac{v(cY+Z)BX}{(A+X)^2}& \displaystyle\frac{vcBX}{A+X} & \displaystyle\frac{vBX}{A+X}-d 
\end{array}\right]\nonumber
\end{equation}

At $F_0$ its eigenvalues are easily found,
$\lambda_1=r$,
$\lambda_2=-s$,
$\lambda_3=-d$.
Since $\lambda_1>0$ the origin is unconditionally unstable. This is a positive result from the conservation point of view, since
the ecosystem will never disappear.

At $F_1$ instead, the characteristic equation factors, to give one explicit eigenvalue
$\lambda_1=-r$, while the remaining ones are the roots of the quadratic
\begin{equation}\label{eqH2}
\lambda^2 + m_1 \lambda + m_0 =0,
\end{equation}
with
$$
m_1=\displaystyle\frac{(s+d)(A+1)-B(wc+v)}{A+1},\quad
m_0=\displaystyle\frac{ds(A+1)-B(sv+cdw)}{A+1} \equiv -W.
$$
We can use Descarte's rule of sign to impose $m_1>0$ and $m_0>0$, so that both roots have negative real part. We thus find, respectively,
$$
A+1>\displaystyle\frac{B(wc+v)}{s+d}, \quad
A+1>\displaystyle\frac{B(sv+cdw)}{ds}.
$$

{\textbf{Remark}}. The feasibility condition for $F_2$ corresponds to $W>0$, so that when $m_0>0$ the only feasible equilibria is $F_1$,
given that $F_0$ is always unstable.

In summary, $F_1$ is locally asymptotically stable if
\begin{equation}\label{eqH4}
A+1>\max\left\{\displaystyle\frac{B(wc+v)}{s+d},\displaystyle\frac{B(sv+cdw)}{ds}\right\}\equiv \displaystyle\frac{B(sv+cdw)}{ds}.
\end{equation}
Note indeed that
$$
\displaystyle\frac{B(sv+cdw)}{ds}>\displaystyle\frac{B(wc+v)}{s+d},
$$
which holds since it reduces to $s^2v+cd^2w>0$, which is true since all parameters are nonnegative.

\begin{figure}[htbp]
\centering
\includegraphics[width=14cm, height= 8cm]{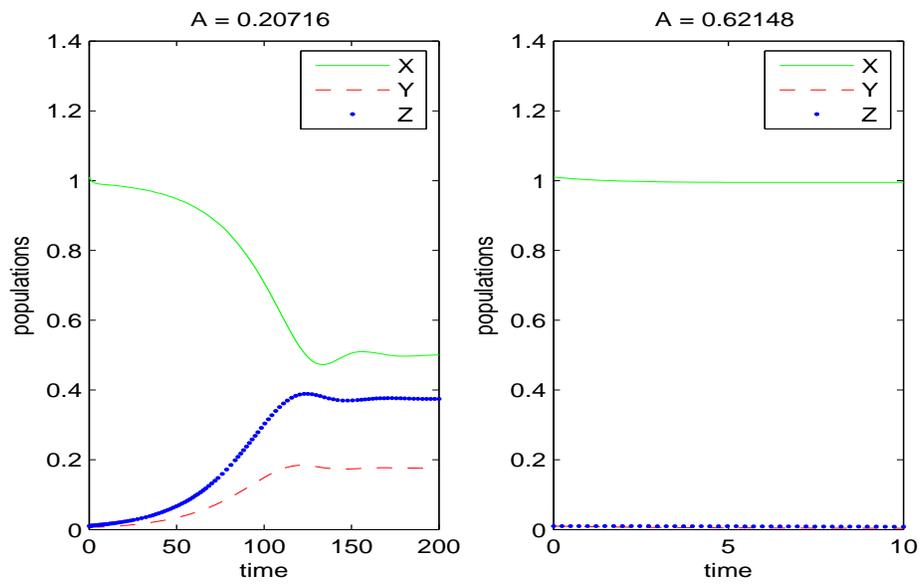}
\caption{Transcritical bifurcation at $F_1$, for the parameter values $r=0.6$, $c=0.38$, $w=0.47$, $s=0.4$, $v=0.5$, $d=0.2$, $B=0.48$,
$A=\displaystyle\frac{i}{2} \times 0.41432$, $i=1,2,3$. On the left the coexistence equilibrium, on the right the equilibrium $E_1$.}
\label{fig5}
\end{figure}

Note that the equilibrium $F_1$ changes stability when the inequality in (\ref{eqH4}) becomes an equality. But this coincides with
the situation that brings $F_2$ to become feasible, see (\ref{sistH1}).
We have thus discovered that there is a transcritical bifurcation, the coexistence equilibrium $F_2$ emanates from the boundary equilibrium $F_1$
when the parameter $B$ attains and crosses the critical value
\begin{equation}\label{transcr}
B^{\dagger}=\frac {ds(A+1)}{sv+cdw}.
\end{equation}
It is illustrated in Figure \ref{fig5}, for the fixed parameter values
$r=0.6$, $c=0.38$, $w=0.47$, $s=0.4$, $v=0.5$, $d=0.2$, $B=0.48$.
The parameter $A$ has then been assigned three different values, namely $A=\displaystyle\frac{i}{2} \, 0.41432$, for $i=1,2,3$.
When $F_1$ is unstable, i.e. for $A=\displaystyle\frac{1}{2} \, 0.41432$, the system settles at the coexistence equilibrium
$(0.5,0.17625,0.375)$.

In summary
\small{ 
\begin{center}
\begin{tabular}{|c|c|c|}
\hline
\small{$0<A+1<\displaystyle\frac{B(wc+v)}{s+d}$} & \small{$\displaystyle\frac{B(wc+v)}{s+d}<A+1<\displaystyle\frac{B(sv+cdw)}{sd}$} & \small{$A+1>\displaystyle\frac{B(sv+cdw)}{ds}$} \\
\hline
$F_1\; unstable$ & $F_1\; unstable$ & $F_1\; stable$\\
\hline
\end{tabular}
\end{center}}
\normalsize

\section{Routh-Hurwitz conditions at coexistence}
To seek for possible interesting behaviors of the system, leading to bifurcations, \cite{KMR},
we need to investigate the eigenvalues of the Jacobian evaluated at the coexistence equilibrium.
For $F_2$, the feasibility condition (\ref{sistH1}) can be recast in the form
\begin{equation}\label{condeq}
0<A\le \displaystyle\frac{V}{ds}.
\end{equation}
The characteristic equation of the Jacobian evaluated at $F_2$ is a cubic
\begin{equation}\label{polch}
\sum_{i=0}^3 a_{3-i} \lambda^i=0,
\end{equation}
with $a_0=1$ and

\begin{eqnarray*}
a_1=\displaystyle\frac{VB(s^2v+cd^2w)+rds[ABQ-W]}{VBQ},\quad a_3=\displaystyle\frac{rdsW}{BQ},\\
a_2=\displaystyle\frac{rds\left\{B(A-1)(s^2v+cd^2w)+(A+1)ds(s+d)+B(wc+v)[W-ds]\right\}}{VBQ}.\\ 
\end{eqnarray*}

We apply the Routh-Hurwitz criterion to (\ref{polch}), imposing
$a_1>0$, $a_3>0$,and
\begin{equation}\label{RH_3}
a_1a_2-a_3>0.
\end{equation}
We now study in terms of the parameter $A$ each one of the first two conditions and the sign of the coefficient $a_2$
to find and exclude intervals for the model parameters arrangements
where the third one possibly does not hold.
The remaining intervals are those in which the stability of $E_2$ may be sought by suitably ``playing" with the parameter $A$.

Observe that
$a_3>0$ is always satisfied when $F_2$ is feasible, in view of the conditions (\ref{condeq}).

\subsection{Study of $a_1$}
For $a_1$ we have the following considerations. The denominator is always strictly positive, in view of (\ref{condeq}).
The numerator is
\begin{equation}\label{numa1}
VB(s^2v+cd^2w)+rds[ds(A+1)+B(sv+cdw)(A-1)].
\end{equation}
The first two factors are always positive, so that the sign depends only on the last term. We study it in terms of the parameter $A$.

If $A \geq 1$ we have easily $a_1>0$. Also, if the bracket is positive, positivity of $a_1$ is once more ensured; this occurs when
\begin{equation}\label{a1pos}
0<\displaystyle\frac{V}{B(sv+cdw)+ds} \leq A <\displaystyle\frac{V}{ds},
\end{equation}
where the inequality on the right is provided by the feasibility condition (\ref{condeq}).

We need still to study the case
$$
0<A<\displaystyle\frac V{B(sv+cdw)+ds}
$$
In this situation, $a_1$ is positive if
\begin{equation}\label{a1posfin}
A>K \equiv \displaystyle \frac V{B(sv+cdw)+ds} \times \displaystyle \frac{rds-B(s^2v+cd^2w)}{rds}.
\end{equation}

\textbf{Remark}. Examining each fraction, we clearly see that $K<1$.
Further, if
$rds>B(s^2v+cd^2w)$ we find $K>0$. Consequently $a_1>0$ holds if (\ref{a1posfin}) is satisfied. Conversely $a_1\leq 0$, if
\begin{equation}\label{a1min}
0<A\leq K.
\end{equation}
In case instead
$rds\leq B(s^2v+cd^2w)$ we have $K\leq 0$ and (\ref{a1posfin}) is always true, so that $a_1>0$.

\textbf{Remark}. Note that the quantity $[B(sv+cdw)-ds](ds)^{-1}$, positive by feasibility of equilibrium
$F_2$, is larger than 1 if and only if $B(sv+cdw)>2ds$.

By combining the considerations for the signs of $a_1$, we have the following four possible situations. In the Table, the interval ranges contain
the possible values of the parameter $A$ and in them we explicitly describe the sign of the coefficient $a_1$,
specifying also the intervals in which the equilibrium $E_2$ is not feasible or does not exist
by the symbol $F_2\,\nexists$.

\begin{enumerate}
\item [(A)]
\begin{tabular}{p{0.8cm}|p{0.8cm}p{0.8cm}|p{0.8cm}p{0.8cm}|p{0.8cm}p{0.8cm}|p{0.8cm}p{0.8cm}|p{0.8cm}}
\multicolumn{2}{c}{K} & \multicolumn{2}{c}{0} & \multicolumn{2}{c}{$\frac{B(sv+cdw)-ds}{B(sv+cdw)+ds}$} &\multicolumn{2}{c}{$\frac{B(sv+cdw)-ds}{ds}$}& \multicolumn{2}{c}{1}\\
\hline
$F_2\,\nexists$&\multicolumn{2}{c|}{$F_2\,\nexists$}&&&&&\multicolumn{2}{c|}{$F_2\,\nexists$}&$F_2\,\nexists$\\
&&&\multicolumn{2}{c|}{$a_1>0$} &\multicolumn{2}{c|}{$a_1>0$} &&&\\
\end{tabular}

\item [(B)]
\begin{tabular}{p{0.8cm}|p{0.8cm}p{0.8cm}|p{0.8cm}p{0.8cm}|p{0.8cm}p{0.8cm}|p{0.8cm}p{0.8cm}|p{0.8cm}}
\multicolumn{2}{c}{K} & \multicolumn{2}{c}{0} & \multicolumn{2}{c}{$\frac{B(sv+cdw)-ds}{B(sv+cdw)+ds}$} &\multicolumn{2}{c}{1}& \multicolumn{2}{c}{$\frac{B(sv+cdw)-ds}{ds}$}\\
\hline
$F_2\,\nexists$&\multicolumn{2}{c|}{$F_2\,\nexists$}&&&&&&&$F_2\,\nexists$\\
&&&\multicolumn{2}{c|}{$a_1>0$} &\multicolumn{2}{c|}{$a_1>0$} &\multicolumn{2}{c|}{$a_1>0$}&\\
\end{tabular}

\item [(C)]

\begin{tabular}{p{0.8cm}|p{0.8cm}p{0.8cm}|p{0.8cm}p{0.8cm}|p{0.8cm}p{0.8cm}|p{0.8cm}p{0.8cm}|p{0.8cm}}
\multicolumn{2}{c}{0} & \multicolumn{2}{c}{K} & \multicolumn{2}{c}{$\frac{B(sv+cdw)-ds}{B(sv+cdw)+ds}$} &\multicolumn{2}{c}{$\frac{B(sv+cdw)-ds}{ds}$}& \multicolumn{2}{c}{1}\\
\hline
$F_2\,\nexists$&&&&&&&\multicolumn{2}{c|}{$F_2\,\nexists$}&$F_2\,\nexists$\\
&\multicolumn{2}{c|}{$a_1<0$}&\multicolumn{2}{c|}{$a_1>0$} &\multicolumn{2}{c|}{$a_1>0$} &&&\\
\end{tabular}

From this Table, observe that $a_1=0$ for $A=K$.

\item [(D)]

\begin{tabular}{p{0.8cm}|p{0.8cm}p{0.8cm}|p{0.8cm}p{0.8cm}|p{0.8cm}p{0.8cm}|p{0.8cm}p{0.8cm}|p{0.8cm}}
\multicolumn{2}{c}{0} & \multicolumn{2}{c}{K} & \multicolumn{2}{c}{$\frac{B(sv+cdw)-ds}{B(sv+cdw)+ds}$} &\multicolumn{2}{c}{1}& \multicolumn{2}{c}{$\frac{B(sv+cdw)-ds}{ds}$}\\
\hline
$F_2\,\nexists$&&&&&&&&&$F_2\,\nexists$\\
&\multicolumn{2}{c|}{$a_1<0$}&\multicolumn{2}{c|}{$a_1>0$} &\multicolumn{2}{c|}{$a_1>0$} &\multicolumn{2}{c|}{$a_1>0$}&\\
\end{tabular}

In this case too, $a_1=0$ for $A=K$.
\end{enumerate} 

\subsection{Study of $a_2$}
To prepare the ground for investigating sufficient conditions leading to the verification of the
last Routh-Hurwitz condition, we begin by studying the sign of $a_2$.

Note that feasibility of $F_2$, $V>0$, implies that the denominator of $a_2$ is positive. For the numerator, we need to analyse the signs of
$A-1$ and of $W-ds=BQ-2ds-dsA$. Requiring them both positive implies clearly that $a_2>0$. This occurs for
\begin{equation}\label{a22}
1\leq A \leq \displaystyle\frac{BQ-2ds}{ds}.
\end{equation}
which is nonempty if and only if $BQ\geq 3ds$.
We need to investigate two cases, corresponding to this last inequality.

\subsubsection{Case 1: $BQ\geq 3ds$.}
As mentioned, if $A\in \left[1, \displaystyle\frac{BQ-2ds}{ds}\right],$ then $a_2>0$.
Otherwise, let
\begin{eqnarray}\label{a23}
M=B(s^2v+cd^2w)+ds[(s+d)-B(wc+v)], \\ \nonumber
H=B(s^2v+cd^2w)+B(wc+v)[2ds-BQ]-(s+d)ds.
\end{eqnarray}
Then, the numerator of $a_2$ becomes $AM-H$ and $a_2>0$ holds if
\begin{equation}\label{a24}
A M>H.
\end{equation}
Observe that $H<M$ strictly, since this inequality explicitly amounts to
\begin{equation}\label{sempreok}
2ds(s+d)+B(wc+v)[BQ-3ds] > 0,
\end{equation}
and the quantity in the last bracket is positive by assumption of Case 1.
Thus, the situations 
$\{M=0,H=0\}$, $\{M=0,H>0\}$, $\{M<0,~H~>~0\}$, $\{M<0,H=0\}$ must all be excluded.
Now, the inequality of Case 1 implies
\begin{equation}\label{ord1}
\displaystyle\frac{BQ-2ds}{ds}\geq 1,
\end{equation}
and furthermore we have
\begin{equation}\label{ord2}
\displaystyle\frac{BQ-2ds}{ds}<\displaystyle\frac{BQ-ds}{ds} \equiv \frac V{ds},
\end{equation}
which is consistent, since the right hand side is positive in view of (\ref{condeq}).

We now analyse the remaining situations.

\textbf{Remark:}
When $M<0$ we always have
\begin{equation}\label{spieg4}
\displaystyle\frac{H}{M}>\displaystyle\frac V{ds},
\end{equation}
since, expanding, we find $Hds<M V$, i.e.
\begin{equation}
B(s^2v+cd^2w)[BQ-2ds ]+ds\{(s+d)BQ-B(wc+v)ds\}>0;\nonumber
\end{equation}
In fact, $ BQ-2ds \geq 0$ since we are in Case 1, namely $BQ\geq 3ds$, and the last brace equals $B(s^2v+cd^2w>0$.
Thus, when $M>0$ we must have the opposite inequality of (\ref{spieg4}), i.e.
\begin{equation}\label{spieg5}
\displaystyle\frac{H}{M}<\displaystyle\frac V{ds}.
\end{equation}

\begin{enumerate}
\item [(1+)]
$
\left\{
\begin{array}{c}
\displaystyle
M>0\\
H>0
\end{array}
\right.
$
\end{enumerate}

Let us define the set
$$
\Omega=\{ ds [B(wc+v)-(s+d)], ds (s+d)+B(wc+v)[BQ-2ds]\}.
$$
Solving the system of inequalities, we have
\begin{eqnarray}
B(s^2v+cd^2w)>\max\,\Omega
=ds (s+d)+B(wc+v)[BQ-2ds]
\end{eqnarray}
so that
\begin{eqnarray}\label{MH}
M=B(s^2v+cd^2w)-\min \Omega\\
H=B(s^2v+cd^2w)-\max \Omega. \nonumber
\end{eqnarray}
Since both are positive, we find $0<\displaystyle\frac{H}{M}<1$.

Thus, there is only one possible arrangement of the various quantities. If $A$ falls in one of the intervals below, the sign of $a_2$ is determined
as in the following Table, since
$a_2>0$ if and only if $A>\displaystyle\frac{H}{M}$.

$$
$$
\begin{tabular}{p{0.8cm}|p{0.8cm}p{0.8cm}|p{0.8cm}p{0.8cm}|p{0.8cm}p{0.8cm}|p{0.8cm}p{0.8cm}|p{0.8cm}}
\multicolumn{2}{c}{0} & \multicolumn{2}{c}{$\frac{H}{M}$} & \multicolumn{2}{c}{1} &\multicolumn{2}{c}{$\frac{B(sv+cdw)-2ds}{ds}$}& \multicolumn{2}{c}{$\frac{B(sv+cdw)-ds}{ds}$}\\
\hline
$F_2\,\nexists$&&&&&&&&&$F_2\,\nexists$\\
&\multicolumn{2}{c|}{$a_2<0$}&\multicolumn{2}{c|}{$a_2>0$} &\multicolumn{2}{c|}{$a_2>0$} &\multicolumn{2}{c|}{$a_2>0$}&
\end{tabular}
$$
$$

Note that $a_2=0$ for $A=\displaystyle\frac{H}{M}$.

\begin{enumerate}
\item [(2+)]
$
\left\{
\begin{array}{c}
\displaystyle
M>0\\
H=0\nonumber
\end{array}
\right.
$ 
\end{enumerate}

The solution is again the inequality (\ref{sempreok}), which always holds. Further, $\displaystyle\frac{H}{M}=0$ so that in this situation we have
$$$$
\begin{tabular}{p{1.1cm}|p{1.1cm}p{1.1cm}|p{1.1cm}p{1.1cm}|p{1.1cm}p{1.1cm}|p{1.1cm}}
\multicolumn{2}{c}{$0=\frac{H}{M}$} & \multicolumn{2}{c}{1} &\multicolumn{2}{c}{$\frac{B(sv+cdw)-2ds}{ds}$}& \multicolumn{2}{c}{$\frac{B(sv+cdw)-ds}{ds}$}\\
\hline
$F_2\,\nexists$&&&&&&&$F_2\,\nexists$\\
&\multicolumn{2}{c|}{$a_2>0$} &\multicolumn{2}{c|}{$a_2>0$} &\multicolumn{2}{c|}{$a_2>0$}&
\end{tabular}
$$$$

\begin{enumerate}
\item [(3+)] 
$
\left\{
\begin{array}{c}
\displaystyle
M<0\\
H<0\nonumber
\end{array}
\right.
$ 
\end{enumerate}

The solution of these inequalities is
\begin{eqnarray}
B(s^2v+cd^2w)<\min\,\{\!\!\!\!\!&ds&\!\!\!\![B(wc+v)-(s+d)];\nonumber\\
&ds&\!\!\!\!\!(s+d)+B(wc+v)[BQ-2ds]\}=\nonumber\\
=&ds&\!\!\!\![B(wc+v)-(s+d)];\nonumber
\end{eqnarray}
(\ref{MH}) again holds, i.e. $H<M$, but both terms are here negative, so that $\displaystyle\frac{H}{M}>1$ follows.

In summary $a_2>0$ when $A<\displaystyle\frac{H}{M}$. But this last quantity exceeds the value for the feasibility of $F_2$.
Thus $a_2>0$ must always hold, namely
$$$$
\begin{tabular}{p{0.8cm}|p{0.8cm}p{0.8cm}|p{0.8cm}p{0.8cm}|p{0.8cm}p{0.8cm}|p{0.8cm}p{0.8cm}|p{0.8cm}}
\multicolumn{2}{c}{0} & \multicolumn{2}{c}{1} &\multicolumn{2}{c}{$\frac{B(sv+cdw)-2ds}{ds}$}& \multicolumn{2}{c}{$\frac{B(sv+cdw)-ds}{ds}$} & \multicolumn{2}{c}{$\frac{H}{M}$}\\
\hline
$F_2\,\nexists$&&&&&&&\multicolumn{2}{c|}{$F_2\,\nexists$}&$F_2\,\nexists$\\
&\multicolumn{2}{c|}{$a_2>0$}&\multicolumn{2}{c|}{$a_2>0$} &\multicolumn{2}{c|}{$a_2>0$} &&&
\end{tabular}
$$$$

\begin{enumerate}
\item  [(4+)]
$
\left\{
\begin{array}{c}
\displaystyle
M\geq 0\\
H<0\nonumber
\end{array}
\right.
$ 
\end{enumerate}

If $M \rightarrow 0^+$, $\displaystyle\frac{H}{M}=-\infty$, and if $M>0$, then $\displaystyle\frac{H}{M}<0$.
In both cases $a_2>0$ strictly, since (\ref{a24}) is easily seen to hold always.

$$$$
\begin{tabular}{p{0.8cm}|p{0.8cm}p{0.8cm}|p{0.8cm}p{0.8cm}|p{0.8cm}p{0.8cm}|p{0.8cm}p{0.8cm}|p{0.8cm}}
\multicolumn{2}{c}{$\frac{H}{M}$}&\multicolumn{2}{c}{0} & \multicolumn{2}{c}{1} &\multicolumn{2}{c}{$\frac{B(sv+cdw)-2ds}{ds}$}& \multicolumn{2}{c}{$\frac{B(sv+cdw)-ds}{ds}$}\\
\hline
$F_2\,\nexists$&\multicolumn{2}{c|}{$F_2\,\nexists$}&&&&&&&$F_2\,\nexists$\\
&&&\multicolumn{2}{c|}{$a_2>0$}&\multicolumn{2}{c|}{$a_2>0$} &\multicolumn{2}{c|}{$a_2>0$} &
\end{tabular}
$$$$

\subsubsection{Case 2: $BQ< 3ds$.}
In this case (\ref{a22}) does not hold. Since also (\ref{sempreok}) does not hold as well, indeed the last bracket in it is now negative,
it is not possible to assess which one among $H$ and $M$ is the larger. We thus need to examine all seven possible configurations for the signs of $H$ and $M$.
Observe further that in this case it follows
\begin{equation}\label{ord0}
\displaystyle\frac{BQ-2ds}{ds}< 1\nonumber
\end{equation}
and more precisely, if $BQ-2ds<0$ we find
\begin{equation}\label{ord3}
\displaystyle\frac{BQ-ds}{ds}< 1,
\end{equation}
while for $BQ-2ds\in (0,ds)$ it follows
\begin{equation}\label{ord4}
\displaystyle\frac{BQ-ds}{ds}> 1.
\end{equation}

\begin{enumerate}
\item [(1-)]
$
\left\{
\begin{array}{c}
\displaystyle
M>0\\
H>0\nonumber
\end{array}
\right.
$ 
\end{enumerate}

We find $\displaystyle\frac{H}{M}\in(0,1)$ if
$$
2ds(s+d)>B(wc+v)[3ds-BQ]\nonumber
$$
and $\displaystyle\frac{H}{M}\geq 1$ for
$$
2ds(s+d)\leq B(wc+v)[3ds-BQ].\nonumber
$$
Here (\ref{spieg5}) still holds, while we can never have
$$
\displaystyle\frac{H}{M}<\displaystyle\frac{BQ-2ds}{ds}<0
$$
but all the other mutual positions of $HM^{-1}$ and $\displaystyle\frac{BQ-2ds}{ds}$ are possible.
In conclusion, we have the following possibilities.

$$$$
\begin{enumerate}
\item [(a)]
\begin{tabular}{p{0.8cm}|p{0.8cm}p{0.8cm}|p{0.8cm}p{0.8cm}|p{0.8cm}p{0.8cm}|p{0.8cm}p{0.8cm}|p{0.8cm}}
\multicolumn{2}{c}{$\frac{B(sv+cdw)-2ds}{ds}$}&\multicolumn{2}{c}{0} & \multicolumn{2}{c}{$\frac{H}{M}$} &\multicolumn{2}{c}{$\frac{B(sv+cdw)-ds}{ds}$}& \multicolumn{2}{c}{1}\\
\hline
$F_2\,\nexists$&\multicolumn{2}{c|}{$F_2\,\nexists$}&&&&&\multicolumn{2}{c|}{$F_2\,\nexists$}&$F_2\,\nexists$\\
&&&\multicolumn{2}{c|}{$a_2<0$}&\multicolumn{2}{c|}{$a_2>0$}&&&
\end{tabular}

For $A=\displaystyle\frac{H}{M}$ we find $a_2=0$.

\item [(b)]
\begin{tabular}{p{0.8cm}|p{0.8cm}p{0.8cm}|p{0.8cm}p{0.8cm}|p{0.8cm}p{0.8cm}|p{0.8cm}p{0.8cm}|p{0.8cm}}
\multicolumn{2}{c}{0}&\multicolumn{2}{c}{$\frac{H}{M}$} & \multicolumn{2}{c}{$\frac{B(sv+cdw)-2ds}{ds}$} &\multicolumn{2}{c}{1}& \multicolumn{2}{c}{$\frac{B(sv+cdw)-ds}{ds}$}\\
\hline
$F_2\,\nexists$&&&&&&&&&$F_2\,\nexists$\\
&\multicolumn{2}{c|}{$a_2<0$}&\multicolumn{2}{c|}{$a_2>0$}&\multicolumn{2}{c|}{$a_2>0$}&\multicolumn{2}{c|}{$a_2>0$}&
\end{tabular}

For $A=\displaystyle\frac{H}{M}$ it follows $a_2=0$.

\item [(c)]
\begin{tabular}{p{0.8cm}|p{0.8cm}p{0.8cm}|p{0.8cm}p{0.8cm}|p{0.8cm}p{0.8cm}|p{0.8cm}p{0.8cm}|p{0.8cm}}
\multicolumn{2}{c}{0}& \multicolumn{2}{c}{$\frac{B(sv+cdw)-2ds}{ds}$}&\multicolumn{2}{c}{$\frac{H}{M}$}  &\multicolumn{2}{c}{1}& \multicolumn{2}{c}{$\frac{B(sv+cdw)-ds}{ds}$}\\
\hline
$F_2\,\nexists$&&&&&&&&&$F_2\,\nexists$\\
&\multicolumn{2}{c|}{$a_2<0$}&\multicolumn{2}{c|}{$a_2<0$}&\multicolumn{2}{c|}{$a_2>0$}&\multicolumn{2}{c|}{$a_2>0$}&
\end{tabular}

For $A=\displaystyle\frac{H}{M}$, again $a_2=0$.

\item [(d)]
\begin{tabular}{p{0.8cm}|p{0.8cm}p{0.8cm}|p{0.8cm}p{0.8cm}|p{0.8cm}p{0.8cm}|p{0.8cm}p{0.8cm}|p{0.8cm}}
\multicolumn{2}{c}{0}& \multicolumn{2}{c}{$\frac{B(sv+cdw)-2ds}{ds}$}  &\multicolumn{2}{c}{1}&\multicolumn{2}{c}{$\frac{H}{M}$}& \multicolumn{2}{c}{$\frac{B(sv+cdw)-ds}{ds}$}\\
\hline
$F_2\,\nexists$&&&&&&&&&$F_2\,\nexists$\\
&\multicolumn{2}{c|}{$a_2<0$}&\multicolumn{2}{c|}{$a_2<0$}&\multicolumn{2}{c|}{$a_2<0$}&\multicolumn{2}{c|}{$a_2>0$}&
\end{tabular}

For $A=\displaystyle\frac{H}{M}$, once  more $a_2=0$.
\end{enumerate}

\begin{enumerate}
\item  [(2-)]
$
\left\{
\begin{array}{c}
\displaystyle
M>0\\
H=0\nonumber
\end{array}
\right.
$ 
\end{enumerate}

Since $\displaystyle\frac{H}{M}=0$, $a_2>0$ holds always, and we have two different possibilities. 

\begin{enumerate}
\item [(a)]
\begin{tabular}{p{1.1cm}|p{1.1cm}p{1.1cm}|p{1.1cm}p{1.1cm}|p{1.1cm}p{1.1cm}|p{1.1cm}}
\multicolumn{2}{c}{$0=\frac{H}{M}$} & \multicolumn{2}{c}{$\frac{B(sv+cdw)-2ds}{ds}$} &\multicolumn{2}{c}{1}& \multicolumn{2}{c}{$\frac{B(sv+cdw)-ds}{ds}$}\\
\hline
$F_2\,\nexists$&&&&&&&$F_2\,\nexists$\\
&\multicolumn{2}{c|}{$a_2>0$} &\multicolumn{2}{c|}{$a_2>0$} &\multicolumn{2}{c|}{$a_2>0$}&
\end{tabular}
$$$$

\item [(b)]
\begin{tabular}{p{1.1cm}|p{1.1cm}p{1.1cm}|p{1.1cm}p{1.1cm}|p{1.1cm}p{1.1cm}|p{1.1cm}}
\multicolumn{2}{c}{$\frac{B(sv+cdw)-2ds}{ds}$} & \multicolumn{2}{c}{$0=\frac{H}{M}$} &\multicolumn{2}{c}{$\frac{B(sv+cdw)-ds}{ds}$}& \multicolumn{2}{c}{1}\\
\hline
$F_2\,\nexists$&\multicolumn{2}{c|}{$F_2\,\nexists$} &&&\multicolumn{2}{c|}{$F_2\,\nexists$}&$F_2\,\nexists$\\
&&&\multicolumn{2}{c|}{$a_2>0$}&&&
\end{tabular}
$$$$
\end{enumerate}

\begin{enumerate}
\item  [(3-)]
$
\left\{
\begin{array}{c}
\displaystyle
M<0\\
H<0\nonumber
\end{array}
\right.
$ 
\end{enumerate}

Here (\ref{spieg4}) holds always and we have three possibilites, recalling (\ref{ord3}) and (\ref{ord4}), since $a_2<0$ if $A<\displaystyle\frac{H}{M}$.

\begin{enumerate}
\item [(a)]
\begin{tabular}{p{0.8cm}|p{0.8cm}p{0.8cm}|p{0.8cm}p{0.8cm}|p{0.8cm}p{0.8cm}|p{0.8cm}p{0.8cm}|p{0.8cm}}
\multicolumn{2}{c}{$\frac{B(sv+cdw)-2ds}{ds}$}& \multicolumn{2}{c}{0}  &\multicolumn{2}{c}{$\frac{B(sv+cdw)-ds}{ds}$}&\multicolumn{2}{c}{$\frac{H}{M}$}& \multicolumn{2}{c}{1}\\
\hline
$F_2\,\nexists$&\multicolumn{2}{c|}{$F_2\,\nexists$}&&&\multicolumn{2}{c|}{$F_2\,\nexists$}&\multicolumn{2}{c|}{$F_2\,\nexists$}&$F_2\,\nexists$\\
&&&\multicolumn{2}{c|}{$a_2>0$}&&&&
\end{tabular}

\item [(b)] 
\begin{tabular}{p{0.8cm}|p{0.8cm}p{0.8cm}|p{0.8cm}p{0.8cm}|p{0.8cm}p{0.8cm}|p{0.8cm}p{0.8cm}|p{0.8cm}}
\multicolumn{2}{c}{0}& \multicolumn{2}{c}{$\frac{B(sv+cdw)-2ds}{ds}$}  &\multicolumn{2}{c}{1}&\multicolumn{2}{c}{$\frac{B(sv+cdw)-ds}{ds}$}& \multicolumn{2}{c}{$\frac{H}{M}$}\\
\hline
$F_2\,\nexists$&&&&&&&\multicolumn{2}{c|}{$F_2\,\nexists$}&$F_2\,\nexists$\\
&\multicolumn{2}{c|}{$a_2>0$}&\multicolumn{2}{c|}{$a_2>0$}&\multicolumn{2}{c|}{$a_2>0$}&&&
\end{tabular}

\item [(c)]
\begin{tabular}{p{0.8cm}|p{0.8cm}p{0.8cm}|p{0.8cm}p{0.8cm}|p{0.8cm}p{0.8cm}|p{0.8cm}p{0.8cm}|p{0.8cm}}
\multicolumn{2}{c}{$\frac{B(sv+cdw)-2ds}{ds}$}& \multicolumn{2}{c}{0}  &\multicolumn{2}{c}{$\frac{B(sv+cdw)-ds}{ds}$}&\multicolumn{2}{c}{1}& \multicolumn{2}{c}{$\frac{H}{M}$}\\
\hline
$F_2\,\nexists$&\multicolumn{2}{c|}{$F_2\,\nexists$}&&&\multicolumn{2}{c|}{$F_2\,\nexists$}&\multicolumn{2}{c|}{$F_2\,\nexists$}&$F_2\,\nexists$\\
&&&\multicolumn{2}{c|}{$a_2>0$}&&&&
\end{tabular}
\end{enumerate}

\begin{enumerate}
\item  [(4-)]
$
\left\{
\begin{array}{c}
\displaystyle
M\geq 0\\
H<0\nonumber
\end{array}
\right.
$ 
\end{enumerate}

If $M\rightarrow 0^+$, clearly $\displaystyle\frac{H}{M}=-\infty$, otherwise this fraction is negative. Only the quantity $\displaystyle\frac{B(sv+cdw)-2ds}{ds}$
can vary, here, namely we find

\begin{enumerate}
\item [(a)]
\begin{tabular}{p{0.8cm}|p{0.8cm}p{0.8cm}|p{0.8cm}p{0.8cm}|p{0.8cm}p{0.8cm}|p{0.8cm}p{0.8cm}|p{0.8cm}}
\multicolumn{2}{c}{$\frac{H}{M}$}& \multicolumn{2}{c}{0}  &\multicolumn{2}{c}{$\frac{B(sv+cdw)-2ds}{ds}$}&\multicolumn{2}{c}{1}& \multicolumn{2}{c}{$\frac{B(sv+cdw)-ds}{ds}$}\\
\hline
$F_2\,\nexists$&\multicolumn{2}{c|}{$F_2\,\nexists$}&&&&&&&$F_2\,\nexists$\\
&&&\multicolumn{2}{c|}{$a_2>0$}&\multicolumn{2}{c|}{$a_2>0$}&\multicolumn{2}{c|}{$a_2>0$}&
\end{tabular}

\item [(b)]
\begin{tabular}{p{0.8cm}|p{0.8cm}p{0.8cm}|p{0.8cm}p{0.8cm}|p{0.8cm}p{0.8cm}|p{0.8cm}p{0.8cm}|p{0.8cm}}
\multicolumn{2}{c}{$\frac{H}{M}$}&\multicolumn{2}{c}{$\frac{B(sv+cdw)-2ds}{ds}$}& \multicolumn{2}{c}{0}& \multicolumn{2}{c}{$\frac{B(sv+cdw)-ds}{ds}$} &\multicolumn{2}{c}{1}\\
\hline
$F_2\,\nexists$&\multicolumn{2}{c|}{$F_2\,\nexists$}&\multicolumn{2}{c|}{$F_2\,\nexists$}&&&\multicolumn{2}{c|}{$F_2\,\nexists$}&$F_2\,\nexists$\\
&&&&&\multicolumn{2}{c|}{$a_2>0$}&&&
\end{tabular}
\end{enumerate}

\begin{enumerate}
\item  [(5-)]
$
\left\{
\begin{array}{c}
\displaystyle
M<0\\
H=0\nonumber
\end{array}
\right.
$ 
\end{enumerate}

Then $\displaystyle\frac{H}{M}=0$. Since $a_2>0$ for $AM>H$, we find here $a_2<0$ always. But then by (\ref{spieg4}) all
other quantites are negative, so that $F_2$ is never feasible.

\begin{enumerate}
\item  [(6-)]
$
\left\{
\begin{array}{c}
\displaystyle
M<0\\
H>0\nonumber
\end{array}
\right.
$ 
\end{enumerate}

Again, $\displaystyle\frac{H}{M}<0$, and (\ref{spieg4}) implies that all
other quantites are negative, so that $F_2$ is never feasible.

\begin{enumerate}
\item  [(7-)]
$
\left\{
\begin{array}{c}
\displaystyle
M=0\\
H>0\nonumber
\end{array}
\right.
$ 
\end{enumerate}

Here $\displaystyle\frac{H}{M}=+\infty$ and $a_2<0$ always, because $AM~>~H$ never holds. Hence stability can never occur.

There are the following alternatives.
\begin{enumerate}
\item [(a)]
\begin{tabular}{p{0.8cm}|p{0.8cm}p{0.8cm}|p{0.8cm}p{0.8cm}|p{0.8cm}p{0.8cm}|p{0.8cm}p{0.8cm}|p{0.8cm}}
\multicolumn{2}{c}{0}&\multicolumn{2}{c}{$\frac{B(sv+cdw)-2ds}{ds}$}& \multicolumn{2}{c}{1}& \multicolumn{2}{c}{$\frac{B(sv+cdw)-ds}{ds}$} &\multicolumn{2}{c}{$\frac{H}{M}$}\\
\hline
$F_2\,\nexists$&&&&&&&\multicolumn{2}{c|}{$F_2\,\nexists$}&$F_2\,\nexists$\\
&\multicolumn{2}{c|}{$a_2<0$}&\multicolumn{2}{c|}{$a_2<0$}&\multicolumn{2}{c|}{$a_2<0$}&&&
\end{tabular}

\item [(b)]
\begin{tabular}{p{0.8cm}|p{0.8cm}p{0.8cm}|p{0.8cm}p{0.8cm}|p{0.8cm}p{0.8cm}|p{0.8cm}p{0.8cm}|p{0.8cm}}
\multicolumn{2}{c}{$\frac{B(sv+cdw)-2ds}{ds}$}&\multicolumn{2}{c}{0}& \multicolumn{2}{c}{$\frac{B(sv+cdw)-ds}{ds}$}& \multicolumn{2}{c}{1} &\multicolumn{2}{c}{$\frac{H}{M}$}\\
\hline
$F_2\,\nexists$&\multicolumn{2}{c|}{$F_2\,\nexists$}&&&\multicolumn{2}{c|}{$F_2\,\nexists$}&\multicolumn{2}{c|}{$F_2\,\nexists$}&$F_2\,\nexists$\\
&&&\multicolumn{2}{c|}{$a_2<0$}&&&&
\end{tabular}
\end{enumerate}

Finally, the case $H=M=0$, is not considered, since $\displaystyle\frac{H}{M}$ is not well defined.

\subsection{Stability of the equilibrium \emph{$F_2$}}

We now combine the previous analyses to assess the situations in which the third Routh-Hurwitz condition holds, namely when (\ref{RH_3}) is satisfied.
Again several cases will arise that are obtained
by suitably merging the previous results. Note that we can merge the two types of Tables rather easily,
since the knots $1$ and $(BQ-ds)(ds)^{-1}$ appear in both of them.
We again study the two cases $BQ\geq 3ds$ and $BQ<3ds$ separately.

\subsubsection{Case 1: $BQ\geq 3ds$.}
\begin{enumerate}
\item  [(1+,2+)]
$
\left\{
\begin{array}{c}
\displaystyle
M>0\\
H\geq 0\nonumber
\end{array}
\right.
$ 
\end{enumerate}
The particular case $H=0$ will be discussed below each Table.
In the present situation, we recall that
\begin{equation}
\displaystyle\frac{BQ-2ds}{ds}\geq 1,\nonumber
\end{equation}
and
\begin{equation}
\displaystyle\frac{BQ-2ds}{ds}<\displaystyle\frac V{ds}.\nonumber
\end{equation}
In what follows we combine only the Tables relative to $a_2$ with those of $a_1$ for which the above relations hold,
i.e. the second and the fourth ones. From the second one, in which $K<0$, we have that $\displaystyle\frac{H}{M}$ can be
larger or smaller than $\displaystyle\frac V{B(sv+cdw)+ds}$. The former alternative holds if and only if
$$
wc+v<\displaystyle\frac{2ds[B(s^2v+cd^2w)-ds(s+d)]}{ds V-[BQ+ds][2ds-BQ]}.
$$
Note that in the fraction the numerator is positive in view of $M>0$, and the denominator is positive as well, since
$2ds-BQ< 0$.

Combining (B) with (1+) and (2+) we obtain the following two Tables, recalling that $a_3>0$ always. At first:

\noindent
\begin{tabular}{p{0.5cm}|p{0.4cm}p{0.4cm}|p{0.4cm}p{0.4cm}|p{0.4cm}p{0.4cm}|p{0.4cm}p{0.4cm}|p{0.4cm}p{0.4cm}|p{0.4cm}p{0.4cm}|p{0.5cm}}
\multicolumn{2}{c}{K} & \multicolumn{2}{c}{0} &\multicolumn{2}{c}{$\frac{B(sv+cdw)-ds}{B(sv+cdw)+ds}$}& \multicolumn{2}{c}{$\frac{H}{M}$}&\multicolumn{2}{c}{1} &\multicolumn{2}{c}{$\frac{B(sv+cdw)-2ds}{ds}$}&\multicolumn{2}{c}{$\frac{B(sv+cdw)-ds}{ds}$}\\
\hline
$F_2\nexists$&\multicolumn{2}{c|}{$F_2\,\nexists$}&&&&&&&&&&$\;\;\;\;\;\;F_2\,\nexists$\\
&&&\multicolumn{2}{c|}{$a_1>0$} &\multicolumn{2}{c|}{$a_1>0$}&\multicolumn{2}{c|}{$a_1>0$}&\multicolumn{2}{c|}{$a_1>0$}&\multicolumn{2}{c|}{$a_1>0$}&\\
&&&\multicolumn{2}{c|}{$a_2<0$} &\multicolumn{2}{c|}{$a_2<0$}&\multicolumn{2}{c|}{$a_2>0$}&\multicolumn{2}{c|}{$a_2>0$}&\multicolumn{2}{c|}{$a_2>0$}&\\
&&&\multicolumn{2}{c|}{$a_3>0$} &\multicolumn{2}{c|}{$a_3>0$}&\multicolumn{2}{c|}{$a_3>0$}&\multicolumn{2}{c|}{$a_3>0$}&\multicolumn{2}{c|}{$a_3>0$}&\\
\end{tabular}
$$$$
Note that for $\displaystyle\frac{H}{M}=0$ this alternative does not exist, since we would find
$\displaystyle\frac V{B(sv+cdw)+ds}~<~0$, impossible because $BQ<ds$ implies the infeasibility of $F_2$.

Now, in this Table and in all the ones that we will consider from now on, we need to find the intervals in which the third Routh-Hurwitz condition
may be satisfied, (\ref{RH_3}). In view of the fact that $a_3>0$ always, we need therefore to identify the intervals in which $a_1$ and $a_2$
have the same sign.
Among those then, one can search whether the condition (\ref{RH_3}) is satisfied.
Evidently, from the above configuration, in this case we find the interval $[\frac HM, \frac V {ds}]$. This is the candidate where
to try for parameter values that will possibly provide a stable $F_2$. The interval is found for this particular arrangements of the knots, and this information
is also relevant. Therefore, we will use the following representation to denote the solution interval for the parameter $A$, within square brackets,
in the corresponding knots arrangements
$$
K<0<\frac V{BQ+ds}<\left[ \frac{H}{M}<1<\frac{BQ-2ds}{ds}<\frac V{ds}\right].
$$
This notation will be used also in what follows, without rewriting explicitly the summarizing Table beforehand.
This arrangement occurs for combining the cases for $a_1$ and $a_2$, i.e. (B,1+,2+). But the above as mentioned is only one of two possible arrangements
in the same situation. The next one is the following one:
\begin{equation}\label{soln_1}
K<0<\left[ \frac{H}{M}<\frac V{BQ+ds}<1<\frac{BQ-2ds}{ds}<\frac V{ds}\right]
\end{equation}

For $H=0$, the table is the same: no matter how $A$ is chosen, all coefficients are always strictly positive.

As long as $A\leq \displaystyle\frac{H}{M}$, the third Routh-Hurwitz condition does not hold, thus $F_2$ in the first interval is unstable.
From the analyses of the Tables, we infer the possibility of a Hopf bifurcation. Although we do not analytically find the bifurcation value
of the parameter, the numerical experiments verify this conjecture.

Next, from (D) and (1+) and (2+) we have the cases corresponding to
(D,1+,2+).
Let us recall that $K<\displaystyle\frac V{BQ+ds}$. Since here $\displaystyle\frac{H}{M}\in(0,1)$
we have the following three possible situations for $K>0$.
$$
\left[ 0<\frac{H}{M}\right] <\left[ K<\frac V{BQ+ds}<1<\frac{BQ-2ds}{ds}<\frac V{ds}\right].
$$

Note that for $\displaystyle\frac{H}{M}=0$, the first interval simply disappears.

$$
\left[ 0<K\right] < \left[ \frac{H}{M}<\frac V{BQ+ds}<1<\frac{BQ-2ds}{ds}<\frac V{ds}\right]
$$

Here the particular case $\displaystyle\frac{H}{M}=0$ cannot hold, since it implies $0<K<\displaystyle\frac{H}{M}=0$.

$$
\left[ 0<K\right] <\frac V{BQ+ds}<\left[ \frac{H}{M}<1<\frac{BQ-2ds}{ds}<\frac V{ds}\right]
$$

For $\displaystyle\frac{H}{M}=0$ the above situation is impossible. 

\begin{enumerate}
\item  [(3+)]
$
\left\{
\begin{array}{c}
\displaystyle
M<0\\
H<0\nonumber
\end{array}
\right.
$ 
\end{enumerate}

Recall that
(\ref{spieg4}) implies that only $K$ influences the dispositions of these points. We thus find the Tables (B) and (D) for $a_1$,
to which we add $\displaystyle\frac{H}{M}$. For (3+,B) we have"
$$
K<\left[ 0<\frac V{BQ+ds}<1<\frac{BQ-2ds}{ds}<\frac V{ds}\right] <\frac{H}{M}
$$

while for (3+,D) instead we find
\begin{equation}\label{soln_2}
0<\left[ K<\frac V{BQ+ds}<1<\frac{BQ-2ds}{ds}<\frac V{ds}\right] <\frac{H}{M}
\end{equation}

Thus, as long as $A\leq K$, the third Routh-Hurwitz condition clearly does not hold.

\begin{enumerate}
\item  [(4+)]
$
\left\{
\begin{array}{c}
\displaystyle
M\geq 0\\
H<0\nonumber
\end{array}
\right.
$ 
\end{enumerate}
There here only two situations, corresponding to $K$ being positive or negative, i.e. respectively to case (D) and (B).
For (4+,B) we have
$$
\frac{H}{M}<K<\left[ 0<\frac V{BQ+ds}<1<\frac{BQ-2ds}{ds}<\frac V{ds}\right]
$$

For (4+,D) we find instead
\begin{equation}\label{soln_3}
\frac{H}{M}<0<\left[ K<\frac V{BQ+ds}<1<\frac{BQ-2ds}{ds}<\frac V{ds}\right]
\end{equation}

Also for the second situation, as long as $A\leq K$, the third Routh-Hurwitz condition does not hold.

\subsubsection{Case 1: $B(sv+cdw)< 3ds$.}
In this case there are many more possibilities.
Let us recall that
$K<\displaystyle\frac V{BQ+ds}$
and that
$$
\displaystyle\frac V{BQ+ds}<\displaystyle\frac{BQ-ds}{ds},\nonumber
$$
always holds, while for the two quantities
$$
\displaystyle\frac V{BQ+ds}, \quad \displaystyle\frac{BQ-2ds}{ds}
$$
one can be larger or smaller than the other one. In all the following cases, the following situations are always true:
\begin{itemize}
\item $a_3>0$ always;
\item $a_2>0$ if

\begin{itemize}
\item[(*)] $A>\displaystyle\frac{H}{M}$, with $M>0$,
\end{itemize}

\item $a_1>0$ if
\begin{itemize}
\item[(*)] $A>K$, with $rds>s^2v+cd^2w$,
\item[(**)] always, with $rds<s^2v+cd^2w$.
\end{itemize}

\end{itemize}

The possible cases are the following ones.
\begin{enumerate}
\item  [(1-)]
$
\left\{
\begin{array}{c}
\displaystyle
M>0\\
H>0\nonumber
\end{array}
\right.
$ 
\end{enumerate}

We now insert the quantities $\displaystyle\frac V{BQ+ds}$ and $K$ in the Tables of the section
relative to $BQ<3ds$; in each situation several subcases will arise, corresponding to different arrangements of the knots.
For the case (a) we have one of the following alternatives when combined with (A), (1-,a,A)
\begin{itemize}
\item $\displaystyle\frac V{BQ+ds}\in \left(0;\displaystyle\frac{H}{M}\right)$, 
\item $\displaystyle\frac V{BQ+ds}\in \left(\displaystyle\frac{H}{M},\displaystyle\frac V{ds}\right),$
\end{itemize}
while $K< \displaystyle\frac V{BQ+ds}$. 
For the case (a) we have one of the following seven alternatives when combined with (A) or (C).

For the case (a) we have four subcases when combined with (A),namely (1-,a,A)
$$
K< \frac{BQ-2ds}{ds}<0<\frac{V}{BQ+ds}<\left[ \frac{H}{M}<\frac{V}{ds}\right]<1
$$

$$
\frac{BQ-2ds}{ds}<K<0<\frac V{BQ+ds}< \left[ \frac{H}{M}<\frac V{ds}\right] <1
$$

$$
K< \frac{BQ-2ds}{ds}<0< \left[ \frac{H}{M}<\frac V{BQ+ds}<\frac V{ds}\right] <1
$$

$$
\frac{BQ-2ds}{ds}<K<0< \left[ \frac{H}{M}<\frac V{BQ+ds}<\frac V{ds}\right] <1
$$

while for (1-,a,C) we find three alternatives
$$
\frac{BQ-2ds}{ds}<\left[ 0<K\right] <\frac V{BQ+ds}< \left[ \frac{H}{M}<\frac V{ds}\right] <1
$$

$$
\frac{BQ-2ds}{ds}<\left[ 0<K\right] < \left[ \frac{H}{M}<\frac V{BQ+ds}<\frac V{ds}\right] <1
$$

$$
\frac{BQ-2ds}{ds}<\left[ 0< \frac{H}{M}\right] <\left[ K<\frac V{BQ+ds}<\frac V{ds}\right] <1
$$

For the Table (b), the alternatives are
\begin{itemize}
\item $\displaystyle\frac V{BQ+ds}\in \left(0;\displaystyle\frac{H}{M}\right)$, or
\item $\displaystyle\frac V{BQ+ds}\in \left(\displaystyle\frac{H}{M},\displaystyle\frac{BQ-2ds}{ds}\right)$, or
\item $\displaystyle\frac V{BQ+ds}\in \left(\displaystyle\frac{BQ-2ds}{ds},1\right);$
\end{itemize}
giving nine alternatives, three with (B) and the remaining ones with (D). For (1-,b,B) we have
$$
K<0<\frac V{BQ+ds}<\left[ \frac{H}{M}<\frac {BQ-2ds}{ds}<1<\frac V{ds}\right]
$$

$$
K<0<\left[ \frac{H}{M}<\frac V{BQ+ds}<\frac {BQ-2ds}{ds}<1<\frac V{ds}\right]
$$

$$
K<0 <\left[ \frac{H}{M}<\frac {BQ-2ds}{ds}<\frac V{BQ+ds}<1<\frac V{ds}\right]
$$

while for (1-,7,b,D) we find
$$
\left[ 0<K\right] <\frac V{BQ+ds}<\left[ \frac{H}{M}<\frac {BQ-2ds}{ds}<1<\frac V{ds}\right]
$$

$$
\left[ 0<K\right] <\left[ \frac{H}{M}<\frac V{BQ+ds}<\frac {BQ-2ds}{ds}<1<\frac V{ds}\right]
$$

$$
\left[ 0< \frac{H}{M}\right] <\left[ K<\frac V{BQ+ds}<\frac {BQ-2ds}{ds}<1<\frac V{ds}\right]
$$

$$
\left[ 0<K\right] <\left[ \frac{H}{M}<\frac {BQ-2ds}{ds}<\frac V{BQ+ds}<1<\frac V{ds}\right]
$$

$$
\left[ 0<\frac{H}{M} \right] <\left[ K<\frac {BQ-2ds}{ds}<\frac V{BQ+ds}<1<\frac V{ds}\right]
$$

$$
\left[ 0<\frac{H}{M} \right] < \frac {BQ-2ds}{ds}<\left[ K<\frac V{BQ+ds}<1<\frac V{ds}\right]
$$

With the Table (c), the alternatives are
\begin{itemize}
\item $\displaystyle\frac V{BQ+ds}\in \left(0;\displaystyle\frac{BQ-2ds}{ds}\right)$, or
\item $\displaystyle\frac V{BQ+ds}\in \left(\displaystyle\frac{BQ-2ds}{ds},\displaystyle\frac{H}{M}\right)$, or
\item $\displaystyle\frac V{BQ+ds}\in \left(\displaystyle\frac{H}{M},1\right)$;
\end{itemize}
giving again nine different cases. With (1-,c,B) we find

$$
K< 0 <\frac V{BQ+ds} < \frac {BQ-2ds}{ds}< \left[ \frac{H}{M} < 1 <\frac V{ds}\right]
$$

$$
K< \left[ 0 < \frac {BQ-2ds}{ds}\right] < \frac V{BQ+ds} <\left[ \frac{H}{M} < 1 <\frac V{ds}\right]
$$

$$
K< 0< \frac {BQ-2ds}{ds} <\left[ \frac{H}{M} < \frac V{BQ+ds} < 1 <\frac V{ds}\right]
$$

$$
\left[ 0< \frac {BQ-2ds}{ds} < K \right] <\left[ \frac{H}{M} < \frac V{BQ+ds} < 1 <\frac V{ds}\right]
$$

$$
\left[ 0< \frac {BQ-2ds}{ds} < \frac{H}{M} \right] <\left[ K < \frac V{BQ+ds} < 1 <\frac V{ds}\right]
$$

while for (1-,c,D) we have
$$
\left[ 0< K \right] <\frac V{BQ+ds} < \frac {BQ-2ds}{ds}< \left[ \frac{H}{M} < 1 <\frac V{ds}\right]
$$

$$
\left[ 0< K\right] < \frac {BQ-2ds}{ds} < \frac V{BQ+ds} <\left[ \frac{H}{M} < 1 <\frac V{ds}\right]
$$

$$
\left[ 0< \frac {BQ-2ds}{ds} < K\right] < \frac V{BQ+ds} <\left[ \frac{H}{M} < 1 <\frac V{ds}\right]
$$

$$
\left[ 0< K \right] < \frac {BQ-2ds}{ds} <\left[ \frac{H}{M} < \frac V{BQ+ds} < 1 <\frac V{ds}\right]
$$

Finally, the Table (d) gives five arrangements, in view of the following alternatives
\begin{itemize}
\item $\displaystyle\frac V{BQ+ds}\in \left(0;\displaystyle\frac{BQ-2ds}{ds}\right)$, or
\item $\displaystyle\frac V{BQ+ds}\in \left(\displaystyle\frac{BQ-2ds}{ds},1\right)$.
\end{itemize}

With (1-,d,B) we have
$$
K< 0 < \frac V{BQ+ds} < \frac {BQ-2ds}{ds} < 1 < \left[ \frac{H}{M} <\frac V{ds}\right]
$$

$$
K< 0  < \frac {BQ-2ds}{ds} < \frac V{BQ+ds} < 1 < \left[ \frac{H}{M} <\frac V{ds}\right]
$$

while for (1-,d,D) we find
$$
\left[ 0< K\right] < \frac V{BQ+ds} < \frac {BQ-2ds}{ds} < 1 < \left[ \frac{H}{M} <\frac V{ds}\right]
$$

$$
\left[ 0< K \right] < \frac {BQ-2ds}{ds} < \frac V{BQ+ds} < 1 < \left[ \frac{H}{M} <\frac V{ds}\right]
$$

$$
\left[ 0 < \frac {BQ-2ds}{ds} <K \right] < \frac V{BQ+ds} <   1   < \left[ \frac{H}{M} <\frac V{ds}\right]
$$

Note that with
$M=0$ and $H>0$ give $\displaystyle\frac{H}{M}=+\infty$, so that in all arrangements we have $a_2<0$ always, case (7-). Thus, as
already remarked, stability is impossible.

\begin{enumerate}
\item  [(2-)]
$
\left\{
\begin{array}{c}
\displaystyle
M>0\\
H=0\nonumber
\end{array}
\right.
$ 
\end{enumerate}

For the Table (a), we have
\begin{itemize}
\item $\displaystyle\frac V{BQ+ds}\in \left(0;\displaystyle\frac{BQ-2ds}{ds}\right)$, or 
\item $\displaystyle\frac V{BQ+ds}\in \left(\displaystyle\frac{BQ-2ds}{ds},1\right),$
\end{itemize}
giving five arrangements including $K$. Note that $\displaystyle\frac{H}{M}=0$. For (2-,a,B) we find

$$
K< \left[ 0 < \frac V{BQ+ds} < \frac {BQ-2ds}{ds} < 1 <\frac V{ds}\right]
$$

$$
K< \left[ 0 < \frac {BQ-2ds}{ds} < \frac V{BQ+ds} < 1 <\frac V{ds}\right]
$$

For (2-,a,D) we have
$$
0< \left[ K < \frac V{BQ+ds} < \frac {BQ-2ds}{ds} < 1 <\frac V{ds}\right]
$$

$$
0< \left[ K < \frac {BQ-2ds}{ds} < \frac V{BQ+ds} < 1 <\frac V{ds}\right]
$$

$$
0< \frac {BQ-2ds}{ds}< \left[ K < \frac V{BQ+ds} < 1 < \frac V{ds}\right]
$$

For the Table (b), there is only one option,
\begin{itemize}
\item $\displaystyle\frac V{BQ+ds}\in \left(0;\displaystyle\frac V{ds}\right)$
\end{itemize}
giving three possibilities for $K$.

For (2-,b,A) we find
$$
K< \frac {BQ-2ds}{ds} < \left[ 0 < \frac V{BQ+ds} < \frac V{ds} \right] < 1
$$

and
$$
\frac {BQ-2ds}{ds} <K< \left[ 0 < \frac V{BQ+ds} <\frac V{ds} \right] < 1
$$

while for (2-,b,C) we have
$$
\frac {BQ-2ds}{ds} <0< \left[ K < \frac V{BQ+ds} <\frac V{ds} \right] < 1
$$

\begin{enumerate}
\item  [(4-)]
$
\left\{
\begin{array}{c}
\displaystyle
M\geq 0\\
H<0\nonumber
\end{array}
\right.
$ 
\end{enumerate}

For the Table (a), here,
\begin{itemize}
\item $\displaystyle\frac V{BQ+ds}\in \left(0;\displaystyle\frac {BQ-2ds}{ds}\right)$, or 
\item $\displaystyle\frac V{BQ+ds}\in \left(\displaystyle\frac{BQ-2ds}{ds},1\right),$
\end{itemize}
and inserting $K$ we have seven cases.
For (4-,a,B) we find
$$
K< \frac HM< \left[ 0< \frac V{BQ+ds} <\frac {BQ-2ds}{ds} <1   <  \frac V{ds} \right]
$$

$$
\frac HM< K< \left[ 0< \frac V{BQ+ds} <\frac {BQ-2ds}{ds} <1   <  \frac V{ds} \right]
$$

$$
K< \frac HM< \left[ 0< \frac {BQ-2ds}{ds} <\frac V{BQ+ds} <1   <  \frac V{ds} \right]
$$

$$
\frac HM< K< \left[ 0< \frac {BQ-2ds}{ds} <\frac V{BQ+ds} <1   <  \frac V{ds} \right]
$$

while for (4-,a,D) we find
$$
\frac HM< 0< \left[ K< \frac V{BQ+ds} <\frac {BQ-2ds}{ds} <1   <  \frac V{ds} \right]
$$

$$
\frac HM< 0< \left[ K< \frac {BQ-2ds}{ds} <\frac V{BQ+ds} <1   <  \frac V{ds} \right]
$$

$$
\frac HM< 0< \frac {BQ-2ds}{ds}< \left[ K <\frac V{BQ+ds} <1   <  \frac V{ds} \right]
$$

For the Table (b), simply
\begin{itemize}
\item $\displaystyle\frac V{BQ+ds}\in \left(0;\displaystyle\frac V{ds}\right),$
\end{itemize}
thus originating four alternatives for $K$.

For (4-,b,A) we have
$$
K< \frac HM < \frac {BQ-2ds}{ds}< \left[ 0 <\frac V{BQ+ds}  <  \frac V{ds} \right] <1
$$

$$
\frac HM < K< \frac {BQ-2ds}{ds}< \left[ 0 <\frac V{BQ+ds}  <  \frac V{ds} \right] <1
$$

$$
\frac HM < \frac {BQ-2ds}{ds}< K< \left[ 0 <\frac V{BQ+ds} < \frac V{ds}  \right] <  1
$$

For (4-,b,C) we have instead
$$
\frac HM < \frac {BQ-2ds}{ds}< 0< \left[ K <\frac V{BQ+ds} < \frac V{ds}   \right] <  1
$$

\begin{enumerate}
\item  [(3-)]
$
\left\{
\begin{array}{c}
\displaystyle
M< 0\\
H\leq 0\nonumber
\end{array}
\right.
$ 
\end{enumerate}

Here we find $\displaystyle\frac{H}{M}>\displaystyle\frac V{ds}$.
For the Tables (a) and (c), we have
\begin{itemize}
\item $\displaystyle\frac V{BQ+ds}\in \left(0;\displaystyle\frac V{ds}\right),$
\end{itemize}
giving three arrangements for each Table, including $K$. 
Here at times we find no solutions, but we include the cases for completeness sake.
The case (3-,a,A) gives 
$$
K< \frac {BQ-2ds}{ds}< 0 <\frac V{BQ+ds}  \frac V{ds} <  \frac HM< 1
$$

$$
\frac {BQ-2ds}{ds}<K< 0 <\frac V{BQ+ds} < \frac V{ds} <  \frac HM< 1
$$

For (3-,a,C) we find
$$
\frac {BQ-2ds}{ds}< [0 < K ]  <\frac V{BQ+ds} < \frac V{ds} <1 <  \frac HM
$$

The case (3-,c,A) gives no solutions. In fact we find:
$$
K< \frac {BQ-2ds}{ds} < 0  <\frac V{BQ+ds} < \frac V{ds} <1 <  \frac HM
$$

$$
\frac {BQ-2ds}{ds}< \left[ K< 0 \right] <\frac V{BQ+ds}  < \frac V{ds} <  \frac HM< 1
$$

For (3-,c,C) we find
$$
\frac {BQ-2ds}{ds}< \left[ 0< K \right] <\frac V{BQ+ds}  < \frac V{ds} <1 <  \frac HM
$$

For the Table (b) instead,
\begin{itemize}
\item $\displaystyle\frac V{BQ+ds}\in \left(0;\displaystyle\frac{BQ-2ds}{ds}\right),$ or
\item $\displaystyle\frac V{BQ+ds}\in \left(\displaystyle\frac{BQ-2ds}{ds},1\right),$
\end{itemize}
so that we have five arrangements including $K$.

No solutions in some cases as well are found, in particular for (3-,b,B):
$$
K< 0 <\frac V{BQ+ds} < \frac {BQ-2ds}{ds}< 1 < \frac V{ds}  <  \frac HM
$$

and

$$
K< 0 < \frac {BQ-2ds}{ds} <\frac V{BQ+ds} < 1 < \frac V{ds}  <  \frac HM.
$$

For (3-,b,D) the same also occurs, but an admissible interval exists:
$$
\left[ 0< K \right] <\frac V{BQ+ds} < \frac {BQ-2ds}{ds}< 1 <    \frac V{ds}  <  \frac HM
$$

$$
\left[ 0< K \right]< \frac {BQ-2ds}{ds} <\frac V{BQ+ds} < 1< \frac V{ds}  <  \frac HM
$$

$$
\left[ 0< \frac {BQ-2ds}{ds} <K \right]< \frac V{BQ+ds} < 1 < \frac V{ds}  <  \frac HM
$$

In this case, if $H=0$ and $M<0$, all other quantities would be negative.

\section{Simulations}

To illustrate the usefulness of the above analysis, for assessing both the stability
of the coexistence as well as for providing a guideline to find possible Hopf bifurcations, \cite{XL},
we provide the results of some numerical simulations.

\begin{figure}[htbp]
\centering
\includegraphics[width=14cm, height= 8cm]{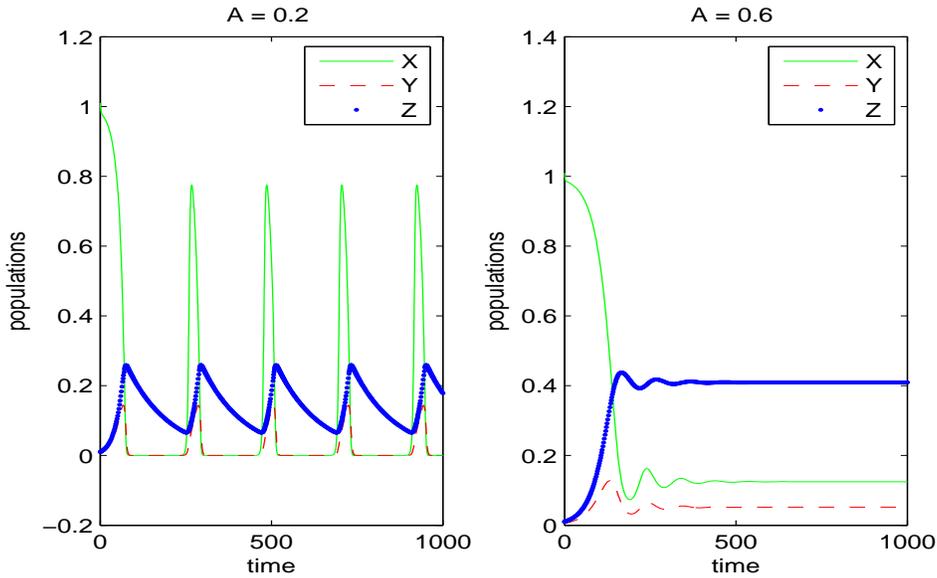}
\caption{Parameter values: $r=0.6$,
$c=0.74$, $w=0.38$, $s=0.48$, $v=0.05$, $d=0.008$, $B=0.85$:
left $A=0.2$; right $A=0.6$.}
\label{fig:H1}
\end{figure}

\textbf{Example 1.} We show at first for instance the results obtained for the situation (\ref{soln_1}).
Figure \ref{fig:H1} reports the system behavior for the parameter values
$r=0.6$, $c=0.74$, $w=0.38$, $s=0.48$, $v=0.05$, $d=0.008$, $B=0.85$.
With this choice, it follows $K\simeq-2.3$,
\begin{eqnarray*}
\displaystyle\frac{H}{M}\simeq 0.36, \quad
\displaystyle\frac{B(sv+cdw)-ds}{B(sv+cdw)+ds}=\frac V{B(sv+cdw)+ds}\simeq 0.71,\\
\displaystyle\frac{B(sv+cdw)-2ds}{ds}\simeq 3.81, \quad
\displaystyle\frac{B(sv+cdw)-ds}{ds}=\frac V{ds}\simeq 4.81.
\end{eqnarray*}
As claimed, we are thus in the situation of (\ref{soln_1}). If we take for $A$ a larger value than the critical value $A=0.4331191029$,
here $A=0.6$, the
coexistence equilibrium is stable, as illustrated in the right plot of Figure \ref{fig:H1}. Taking instead
$A$ in the first interval of (\ref{soln_1}), say $A=0.2$ we see that limit cycles appear, left plot.

\begin{figure}[htbp]
\centering
\includegraphics[width=14cm, height= 8cm]{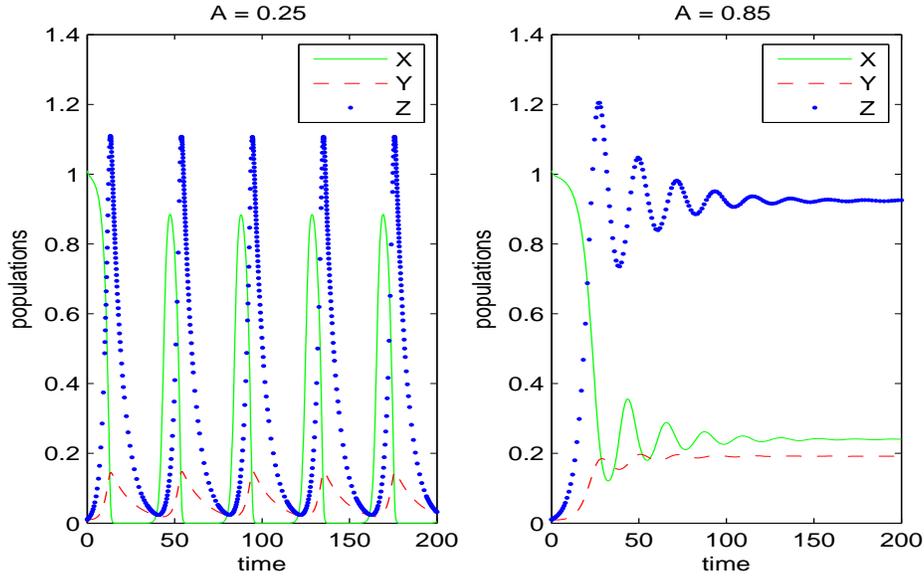}
\caption{Parameter values: $r=0.95$,
$c=0.066$, $w=0.083$, $s=0.075$, $v=0.8$, $d=0.15$, $B=0.84$:
left $A=0.25$; right $A=0.85$.}
\label{fig:H2}
\end{figure}

\textbf{Example 2.} We illustrate now the case of (\ref{soln_2}).
Taking the following parameter values
$r=0.95$,
$c=0.066$, $w=0.083$, $s=0.075$, $v=0.8$, $d=0.15$, $B=0.84$,
we find that $K\simeq 0.41$,
\begin{eqnarray*}
\displaystyle\frac{H}{M}\simeq 15.03, \quad \displaystyle\frac{B(sv+cdw)-ds}{B(sv+cdw)+ds}=\frac V{B(sv+cdw)+ds}\simeq 0.64,\\
\displaystyle\frac{B(sv+cdw)-2ds}{ds}\simeq 2.54, \quad \displaystyle\frac{B(sv+cdw)-ds}{ds}=\frac V{ds}\simeq 3.54.
\end{eqnarray*}
As long as $A\leq K$, here we took $A=0.25$, left plot of Figure \ref{fig:H2}, we find limit cycles.
Past the critical value $A=.6376318460$, the coexistence equilibrium is stable. This is shown on the right plot for $A=0.85$.

\begin{figure}[htbp]
\centering
\includegraphics[width=14cm, height= 8cm]{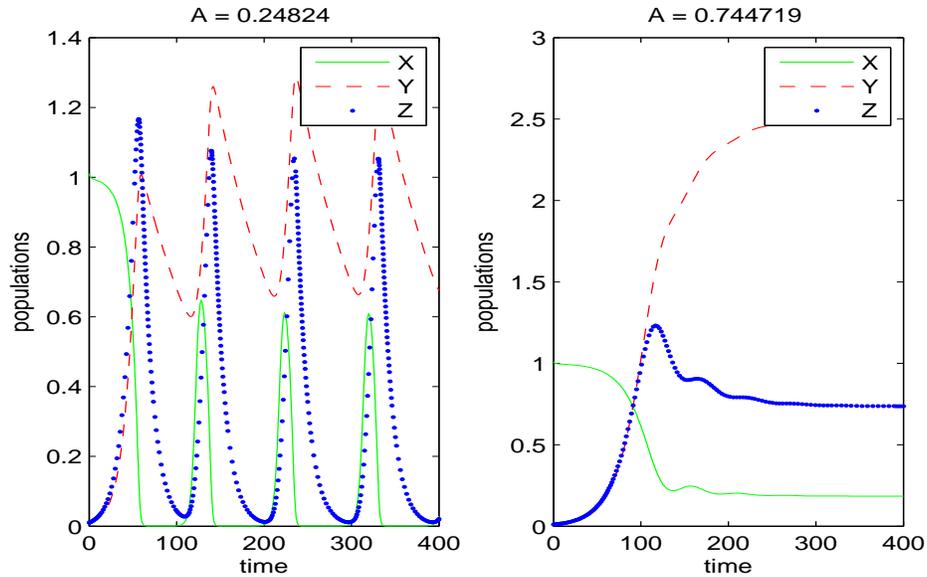}
\caption{Parameter values: $r=0.56$,
$c=0.44$, $w=0.3$, $s=0.01$, $v=0.7$, $d=0.08$, $B=0.23$:
left $A=\frac 12 \times 0.4964791610$; right $A=\frac 32 \times 0.4964791610$.}
\label{fig:H3}
\end{figure}

\textbf{Example 3.} One more instance is shown for the case (\ref{soln_3}).
We take $r=0.56$,
$c=0.44$, $w=0.3$, $s=0.01$, $v=0.7$, $d=0.08$, $B=0.23$.
This choice gives $K\simeq 0.24$,
\begin{eqnarray*}
\displaystyle\frac{H}{M}\simeq -2.54, \quad \displaystyle\frac{B(sv+cdw)-ds}{B(sv+cdw)+ds}=\frac V{B(sv+cdw)+ds}\simeq 0.67,\\
\displaystyle\frac{B(sv+cdw)-2ds}{ds}\simeq 3.05, \quad \displaystyle\frac{B(sv+cdw)-ds}{ds}=\frac V{ds}\simeq 4.05.
\end{eqnarray*}
Now for values of $A$ below the threshold $A=0.4964791610$, here we take the half of that value, sustained oscillations arise,
while for larger values, we take one and a half that critical value, the coexistence equilibrium is stable. These results are shown in Figure \ref{fig:H3}.

\newpage
\section{Conclusions}

From the conservationist point of view,
a nice feature of the ecosystem presented here is that it can never disappear, as the origin is always unstable.
Furthermore, when the prey-only equilibrium is unstable, the system is permanent, \cite{BFW}.

In this system with a response function that models the feeding satiation, oscillations have been shown to arise, through
an in depth investigation of the possible signs of the coefficients in the characteristic equation related to the coexistence
equilibrium. Clearly, the full ecosystem can thrive also at a stable steady state. The result on limit cycles parallels the
one found for the corresponding situation in which rather it is the prey that are genetically distinct, \cite{EV12}.
The model in which genetic differences in predators combine with a standard quadratic response function instead does not
show this feature, \cite{VV}. The models with different genotypes in the predators further show that the coexistence
equilibrium emanates from the prey-only equilibrium under specific system's features, (\ref{transcr}), due to the presence of a
transcritical bifurcation.

Another interesting feature common to these models,
is that it is not possible to have equilibrium with just one genotype.
At first this result is quite surprising, but its more careful analysis
shows that it is inherent in the model assumptions. In fact, new genotypes
can arise from an original genotype. This fact is modeled in the reproduction
terms of the system, compare the last two equations of (\ref{H1}). In fact both
$y$ and $z$ populations have offsprings also belonging to the other population.
Even if one of them gets extinguished at some instant in time, it will be eventually
replenished by the mutations occurring in the other one.
The critical value of the parameter $B$ in condition (\ref{H1}) acts also as an
indicator of the predators invasion of the system. This result is in line with similar ones that hold for the two
models presented in \cite{EV12,VV}.

The conclusion of \cite{VV} that
genetical diversity of the population may affect in a different way the ecosystem, depending on which trophic level it lies,
appears here however more tied to the way the response function that is assumed to hold in the system.

\vspace{10pt} \noindent
{\bf Acknowledgements:} The authors thank Professor G. Badino (Dipartimento di Scienze della Vita e Biologia dei Sistemi, Univ. of Torino)
for a very useful discussion on the matters of this research.
This research was partially supported by
the project ``Metodi numerici in teoria delle popolazioni''
of the Dipartimento di Matematica ``Giuseppe Peano''.


\begin{thebibliography}{11}



\bibitem{BC} E.~Beltrami, T.~O.~Carroll, Modelling the role of viral disease
in recurrent phytoplankton blooms, {\em J. Math. Biol.}~32, 1994, 857--863.

\bibitem{BCC} F. Brauer, C. Castillo-Chavez, \emph{Mathematical Models in Population Biology and Epidemiology}, Springer, 2002.

\bibitem {BFW} G.~Butler, I.~Freedman, P.~Waltman, Uniformly persistent systems, {\em Proc. Amer.
Math. Soc.}~96, 1986, 425--430.

\bibitem{CA} Chattopadhyay, J., Arino, O., 1999. A predator–prey model with disease in the prey. Nonlinear Analysis
36, 747--766.

\bibitem{CB} Chattopadhyay, J., Bairagi, N., 2001. Pelicans at risk in Salton Sea – an eco–epidemiological study.
Ecological Modelling 136, 103--112.

\bibitem{CP} Chattopadhyay, J., Pal, S., 2002. Viral infection on phytoplankton–zooplankton system–a mathematical
model. Ecological Modelling 151, 15--28.

\bibitem {Cal} E.~J.~Calabrese, {\em Ecogenetics: genetic variation in susceptibility to environmental agents},
Environmental Science and Technology, New York, 1984.

\bibitem {Cu} M.~R.~Cummings, {\em Human Heredity: Principles and Issues}, Thomson Brooks/Cole, Belmont, 2000.

\bibitem{LC07} S. Gakkhar, B. Singh, R. K. Naji, \emph{Dynamical behavior of two predators
competing over a single prey}, 
Biosystems 90, 808--817, 2007.

\bibitem{HF}  K.~P.~Hadeler, H.~I.~Freedman, Predator-prey
populations with parasitic infection, {\em J. of Math. Biology}~27, 1989, 609--631.

\bibitem{LC06} V. K$\check{\textrm{r}}$ivan, J. Eisner, \emph{The effect of the Holling type II functional response on apparent competition},
Theoretical Population Biology 70, 421--430, 2006.

\bibitem{KMR} Yu. A. Kuznetsov, S. Muratori, S. Rinaldi, \emph{Bifurcation and chaos in a periodic predator-prey model},
International Journal of Bifurcation and Chaos 2, 117-128, 1992.

\bibitem{LC03} X. Liu, L. Chen, \emph{Complex dynamics of Holling type II
Lotka-Volterra predator-prey system with impulsive perturbations on the predator}, Chaos, Solitons and Fractals
6, 311--320, 2003.

\bibitem{MPV} H.M~alchow, S.~Petrovskii, E.~Venturino, {\em Spatiotemporal patterns in Ecology and Epidemiology},
CRC, Boca Raton, 2008.

\bibitem{Mu} J.D. Murray, \emph{Mathematical Biology. An Introduction. Third Edition}, Springer-Verlag, 2002.

\bibitem{EV92} E.~Venturino, The influence of diseases on Lotka-Volterra systems, {\em IMA preprint \#951}, Minneapolis, MN, 1992.

\bibitem{EV94} E.~Venturino, The influence of diseases on Lotka-Volterra systems, {\em Rocky
Mountain Journal of Mathematics}~24, 1994, 381--402.

\bibitem{EV95} E.~Venturino, Epidemics in predator-prey models: disease among the prey, in O.
Arino, D. Axelrod, M. Kimmel, M. Langlais: {\em Mathematical Population
Dynamics: Analysis of Heterogeneity, Vol. one: Theory of Epidemics},
Wuertz Publishing Ltd, Winnipeg, Canada, p. 381-393, 1995.

\bibitem{EV01}  E. Venturino, The effects of diseases on competing species, {\em Math. Biosc.}, v.
174, p. 111-131, 2001.

\bibitem{EV02}  E. Venturino, Epidemics in predator-prey models: disease in the predators,
{\em IMA Journal of Mathematics Applied in Medicine and Biology} 19,
185-205, 2002.

\bibitem{EV07}  E. Venturino, How diseases affect symbiotic communities, Math. Biosc. 206, 11-30, 2007.

\bibitem{EV12}
E.~Venturino, An ecogenetic model, {\em Appl. Math. Letters}~25, 2012, pp.~1230--1233.

\bibitem{VV} C. Viberti, E. Venturino, A Predator-Prey Model with Genetically Distinguishable Predators, in A. Kanarachos, N. E. Mastorakis (Editors)
Recent Advances in Environmental Sciences,
Proceedings of the 9th International Conference on Energy, Environment, Ecosystems and Sustainable Development (EEESD'13), Lemesos, Cyprus, March
21st-23rd 2013, WSEAS Press, p. 87-92, ISSN 2227-4359, ISBN 978-1-61804-167-8.

\bibitem{XL} D. Xiao, W. Li, \emph{Limit Cycles for the Competitive Three Dimensional Lotka-Volterra System}, Journal of Differential Equations, 2000.

\end{thebibliography}
\end{document}